\documentclass[12pt,twoside]{article}
\usepackage[all]{xy}
\usepackage{a4,amsmath,amssymb,amsfonts,amscd,mathrsfs}
\reversemarginpar
\newcommand{\kkk}[1]{}
\newcounter{Abschnitt}[section]
\newcommand{\neu}[1]{\protect\refstepcounter{Abschnitt}\protect
   \label{#1}\vspace{1ex}
   {\bf \protect\arabic{section}.\protect\arabic{Abschnitt}}
                     \kkk{#1}}

\newcommand{\Db}{\mathrm{D^b}}
\newcommand{\DD}{\mathrm{D}}

\newcommand{\Ext}{{\rm Ext}}
\newcommand{\FM}{{\rm FM}}

\newcommand{\Hom}{{\rm Hom}}
\newcommand{\Jac}{{\rm Jac}}

\newcommand{\NS}{{\rm NS}}
\newcommand{\Pic}{{\rm Pic}}

\newcommand{\Quot}{{\rm Quot}}

\newcommand{\SU}{{\rm SU}}

\newcommand{\U}{{\rm U}}

\newcommand{\im}{{\rm im}}

\newcommand{\rk}{{\rm rk}}

\newcommand{\Ecal}{{\cal E}}

\newcommand{\Jcal}{{\cal J}}

\newcommand{\Ocal}{{\cal O}}
\newcommand{\Pcal}{{\cal P}}

\newcommand{\ndop}{{\mathbb N}}

\newcommand{\sext}{{{\cal E} xt}}

\newcommand{\dual}{^\lor}

\newcommand{\Rar}{\Rightarrow}
\newcommand{\rarpa}[1]{\stackrel{#1}{\rightarrow}}

\newcommand{\Lrar}{\Leftrightarrow}

\newcommand{\proof}{{\bf Proof: }}

\newcommand{\qed}{{ \hfill $\Box$}}

\newcommand{\tip}{{\vspace{0.5em} }}

\author{Georg Hein}
\date{July 10, 2006}
\begin{document}
\title{Raynaud's vector bundles and base points of the generalized Theta divisor}
\maketitle

\begin{abstract}
We study base points of the generalized $\Theta$-divisor on the
moduli space of vector bundles on a smooth algebraic curve $X$
of genus $g$ defined over an algebraically closed field. 
To do so, we use the derived categories $\Db(\Pic^0(X))$ and $\Db(\Jac(X))$
and the equivalence between them given by the Fourier-Mukai transform
$\FM_\Pcal$ coming from the Poincar\'e bundle.
The vector bundles $P_m$ on the curve $X$
defined by Raynaud play a central role in this description.
Indeed, we show that $E$ is a base point of the generalized
$\Theta$-divisor, if and only if there exists a nontrivial homomorphism
$P_{\rk(E)g+1} \to E$.
\end{abstract}

\section{Introduction}
Let $X$ be a smooth projective curve of genus $g \geq 1$
defined over an algebraically closed field of arbitrary characteristic.
If $E$ is a sheaf on $X$ such that $h^0(X,E \otimes L)=0=h^1(X,E \otimes
L)$ for a line bundle $L$ of degree zero, then it is an easy exercise to
check that $E$ is a semistable vector bundle with slope
$\mu(E):=\frac{\deg(E)}{\rk(E)}=g-1$. However, not for all semistable
vector bundles $E$ of that slope does such a line bundle $L$ exist.
Those stable vector bundles form the base locus of the generalized
$\Theta$-divisor. Therefore, we say that a vector bundle $E$ has no
$\Theta$-divisor if $H^0(E \otimes L) \ne 0$ for all line bundles $L$ of
degree zero.

The first examples of such vector bundles were constructed by Raynaud in
\cite{Ray} giving one of the first applications of Fourier-Mukai
transforms developed in \cite{Muk}.
He constructed a sequence of semistable vector bundles $P_m$ for $m \in \ndop$
having the property that $H^0(L \otimes P_m)\ne 0$ for all line bundles $L$ of
degree zero.

A vector bundle $E$ with a nontrivial morphism $P_m \to E$ yields an
example of a vector bundle without $\Theta$-divisor.
Our aim is to show these are all examples.
We show in theorem \ref{PM}.\ref{pm-prop} that a vector bundle $E$ has
no $\Theta$-divisor if and only if $\Hom(P_m,E) \ne 0$ for $m \gg0$.
In Theorem \ref{RANK}.\ref{rank-thm2} we show that the condition on $m$
can be made effective.
In Section \ref{BAD} we define minimal vector bundles to be the minimal
ones having no $\Theta$-divisor and show in Theorem
\ref{BAD}.\ref{bad-thm} how they are related to ample divisors on the
Picard group $\Pic^0(X)$.
In section \ref{SPEC} we use the theory of spectral curves to show that
we obtain Raynaud bundles $P_{m,R}$ having the property that for any
vector bundle $E$ on $X$ the conditions
$H^0(E \otimes F) \ne 0$ for all vector bundles $F$ of rank $R$ and
degree zero is equivalent to $\Hom(P_{m,R},E) \ne 0$. For more details
see Theorem \ref{SPEC}.\ref{spec-thm}. Finally,
in the last section we provide some more applications.\\

{\bf Acknowledgment.} I would like to thank K.~Altmann and D.~Ploog for
excellent discussions at their interesting university.

\section{The vector bundles $P_m$}\label{PM}\kkk{PM}
\neu{pm-notation}{\bf Notations.\/}
Let $\Pcal_X$ be a Poincar\'e bundle on $X \times \Pic^{0}(X)$.
This line bundle is the pull back of a Poincar\'e bundle $\Pcal$ on $\Jac(X)
\times \Pic^{0}(X)$.
Here $\Jac(X)$ denotes the dual abelian variety of $\Pic^{0}(X)$ which
is isomorphic
to to the principally polarized abelian variety $\Pic^{0}(X)$.
Let $x_0$ be a geometric point on $X$. We denote its image in $\Jac(X)$
also by $x_0$.
With this point as zero $\Jac(X)$ becomes a group.
We denote the map assigning an element its inverse by $[-1]_{\Jac(X)}$.
The Poincar\'e bundle $\Pcal$ becomes unique, when we require
$\Pcal|_{\{ x_0 \} \times \Pic^{0}(X)} \cong \Ocal_{\{ x_0 \} \times
\Pic^{0}(X)}$.
On $\Pic^{0}(X)$ we have the Theta divisor
$\Theta_\Pic := \{ [ L] \in \Pic^{0}(X) \, |
\, h^1(X, L \otimes L_{g-1}) >0 \}$ where $L_{g-1}:=\Ocal_X((g-1)x_0)$
is a fixed line bundle of degree $g-1$ on $X$. 
Resuming, we have the following commutative diagram of varieties and
morphisms:
\[\xymatrix{
X \times \Pic^{0}(X) \ar[d]_{q_X} \ar[r]^-{\iota_\Pic}
\ar@(ur,ul)[rr]^{p_X}
& \Jac(X) \times \Pic^{0}(X) \ar[r]^-{p} \ar[d]_-{q} & \Pic^{0}(X)\\
X \ar@{^(->}[r]^-{\iota} & \Jac(X)\\
}\]
We will use the following line bundles
\begin{center}
\begin{tabular}{ll}
$\Pcal \in \Pic(\Jac(X) \times \Pic^{0}(X)) $ & the Poincar\'e bundle\\
%$\Pcal\dual \in \Pic(\Jac(X) \times \Pic^{0}(X) $ & its dual \\
$\Pcal_X = \iota_\Pic^* \Pcal $ & the universal bundle on $X \times
\Pic^{0}(X)$ \\
%$\Pcal\dual_X $ & its dual\\
$\Ocal_{\Pic}(\Theta_\Pic) $ & the Theta line bundle on the Picard group\\
\end{tabular}
\end{center}

\neu{pm-pm}{\bf Definition of the vector bundle $P_m$.\/}
Let $m$ be a positive integer.
We define the vector bundle $P_m$ to be the following direct image sheaf:
\[ P_m := \iota^* [-1]^*_{\Jac(X)} R^gq_* (\Pcal \otimes p^* \Ocal_\Pic(-m \cdot
\Theta_\Pic) )\,.\]
From Kodaira vanishing theorem and Serre duality we conclude the
vanishing of the
direct image sheaves
$R^iq_* (\Pcal \otimes p^* \Ocal_\Pic(-m \cdot \Theta_\Pic) )$
for $i < g$. Therefore, we can use the Riemann-Roch theorem to calculate
the
numerical invariants of $P_m$.

\neu{pm-num-inv}{\bf Lemma. \/}{\em
The vector bundle $P_m$ is semistable with numerical invariants
\[ \rk(P_m)= m^g \qquad \deg(P_m) = g \cdot  m^{g-1} \qquad \mu(P_m)
= \frac{g}{m}\,.\]
}

\proof This is a direct consequence of 3.1.~in Raynaud's article
\cite{Ray}.
The only difference is the application of $[-1]^*_{\Jac(X)}$ in our
definition of $P_m$. Since $C$ and $[-1]^*_{\Jac(X)}C$ are 0-homologous
cycles this morphism does
not change the numerical invariants of the bundles.
\qed

\neu{pm-chain}{\bf Lemma. \/}{\em
For any $m>0$ there exists a canonical surjection
$\xymatrix{P_{m+1} \ar@{->>}[r] & P_m}$.}

\tip
\proof
Since $(\Pic^{0}(X),\Theta_\Pic)$ is principally polarized,
there exist (up to a scalar multiple) a unique morphism 
$\Ocal_{\Pic^{0}(X)}(-(m+1)\Theta_\Pic)  \rarpa{\psi}
\Ocal_{\Pic^{0}(X)}(-m\Theta_\Pic)$.
The morphism $\psi$ is injective and its cokernel has support on the
divisor $\Theta$.
Thus, we obtain the canonical surjection on $\Jac(X)$
\[\xymatrix{R^gq_{*} (\Pcal \otimes p^* \Ocal_\Pic(-(m+1) \cdot
\Theta_\Pic) ) \ar@{->>}[r] &
R^gq_{*} (\Pcal \otimes p^* \Ocal_\Pic(-m \cdot \Theta_\Pic) )
\,.} \]
Pulling this surjection via $\iota^* [-1]^*_{\Jac(X)}$
back to $X$ yields the asserted surjection.
\qed

\neu{pm-prop}{\bf Theorem.\/} {\em
For a coherent sheaf $F$ on $X$ we have the equivalence:
\[( \mbox{For all } [L] \in \Pic^0(X) \mbox{ we have } \Hom(L,F) \ne 0 )
\Lrar ( \Hom(P_m,F) \not= 0
\mbox{ for all } m\gg 0) \,.\] }
\proof
Suppose $\Hom(L,F)=0$ for a line bundle $[L] \in \Pic^0(X)$.
By semicontinuity this holds for all $[L] \in U \subset \Pic^0(X)$
for some open set $U$.
Thus, the sheaf
$p_{X*}(\Pcal_X \otimes q_X^* F)$ is torsion. Since the direct image
sheaf
$p_{X*}(\Pcal_X \otimes q_X F)$
is torsion free, we conclude that $p_{X*}(\Pcal_X \otimes q_X F)=0$.
Expressing this in terms of the Fourier-Mukai transform
$\FM_\Pcal:\Db(\Jac(X)) \to \Db(\Pic^{0}(X))$
we obtain the equality $\FM_\Pcal(\iota_*F) = G[-1]$ for a sheaf
$G$ on $\Pic^{0}(X)$.
Consequently, we have
\[\Hom_{\Db(\Pic^{0}(X))}(\Ocal_{\Pic^{0}(X)}(-m \Theta),\FM_\Pcal(\iota_*F))=0\]
for all integers $m$.
The Fourier transform $\FM_{\Pcal[g]} : \Db(\Pic^{0}(X)) \to \Db(\Jac(X))$
given by the kernel $\Pcal[g]$ yields
\[\Hom_{\Db(\Jac(X))}(\FM_{\Pcal[g]}( \Ocal_{\Pic^{0}(X)}(-m
\Theta)),[-1]^*\iota_*F)= 0\,.\]
Thus, we have $\Hom_{\Db(X)}(\iota^*[-1]^* \FM_{\Pcal[g]}(
\Ocal_{\Pic^{0}(X)}(-m \Theta)), F)=0$.
By the definition of $P_m$ this equals $\Hom_{\Ocal_X}(P_m,F)=0$ for all
integers $m \geq 1$.

If $\Hom(L,F) \ne 0$ for all $[L] \in \Pic^0(X)$, 
then, we have $F^0:=p_{X*}(\Pcal_X \otimes q_XF)$
is a nontrivial coherent sheaf on $\Pic^{0}(X)$.
Thus, for $m\gg 0$ we have $F^0(m \Theta)$ has global sections, or
equivalently
$\Hom(\Ocal_{\Pic^{0}(X)}(-m \Theta),F^0) \ne 0$.
We conclude
\[ \Hom_{\Db(\Pic^{0}(X))}(\Ocal_{\Pic^{0}(X)}(-n
\Theta),\FM_\Pcal(\iota_*F))\ne 0\,.\]
As before, this yields $\Hom_{\Ocal_X}(P_m,F)\ne 0$ for $m \gg 0$.
\qed

\section{Semistable vector bundles without
$\Theta$-divisor}\label{RANK}\kkk{RANK}

\neu{rank-pre}{\bf Preliminaries.}
From \ref{RANK}.\ref{rank-pre} to Corollary
\ref{RANK}.\ref{rank-h0} in this section $E$ denotes a vector bundle of rank
$r$ and the property that $\chi(F) \leq 0$ for all subsheaves $F \subset E$.
The last condition is equivalent to $\mu_{\rm max}(E) \leq g-1$.
Furthermore, we assume that $E$ has no $\Theta$-divisor.
Thus, the $\Ocal_{\Pic^{0}(X)}$-sheaves $\FM_\Pcal^0(E):=
R^0p_{X*}(\Pcal \otimes q_X^*E)$ and 
$\FM_\Pcal^1(E):= R^1p_{X*}(\Pcal \otimes q_X^*E)$ have positive rank.
Let $D=P_1+P_2+ \ldots + P_{g-1}$ be a reduced divisor of degree $g-1$ on
$X$.
Applying the functor $R^\bullet p_{X*}(\Pcal \otimes q_X^* 
\underline{\quad}\, )$
to the short exact sequence $0 \to E(-D) \to E \to \Ocal_D^{\oplus r} \to
0$
gives a long exact sequence
\[\FM_\Pcal^0(E(-D)) \to \FM_\Pcal^0(E) \to \FM_\Pcal^0(\Ocal_D^{\oplus
r}) \to 
\FM_\Pcal^1(E(-D)) \to \FM_\Pcal^1(E) \to \FM_\Pcal^1(\Ocal_D^{\oplus r})
\]
of coherent sheaves on $\Pic^{0}(X)$.

\neu{rank-monad}{\bf Lemma. \/}{\em The sheaves $\FM_\Pcal^0(E(-D))$ and
$\FM_\Pcal^1(\Ocal_D^{\oplus r})$ are zero. The direct image sheaf
$\FM_\Pcal^0(\Ocal_D^{\oplus r})$ is isomorphic to a direct sum of
$r(g-1)$ line bundles
numerically equivalent to $\Ocal_{\Pic^{0}(X)}$.
}

\proof
Let $L$ be a line bundle of degree $0$.
If $h^0(E(-D) \otimes L)>0$, then we have a non trivial morphisms 
$L\dual \to E(-D)$. Since $\mu_{\rm max}(E(-D)) \leq 0$
this implies that $\mu_{\rm max}(E(-D))=0$ and that $L\dual$ is 
is a direct summand of the graduated object associated to the subbundle
$E_{\rm max} \subset E(-D))$ of maximal slope.
We deduce, that there exist at most $r$ line bundles $L \in \Pic^{0}(X)$
with $h^0(E(-D) \otimes L)>0$. Consequently, the torsion free sheaf
$\FM_\Pcal^0(E(-D))$ is zero.

Since the support of $\Ocal_D$ is zero dimensional, there
are no higher direct images. In particular, $\FM_\Pcal^1(\Ocal_D^{\oplus
r})=0$.
$D$ is a reduced divisor, thus $\Ocal_D = \oplus_{i=1}^{g-1} k(P_i)$.
Thus,
$\FM^0_\Pcal(\Ocal_D^{\oplus r}) = \left( \oplus_{i=1}^{g-1} \FM^0_\Pcal(
k(P_i) ) \right)^{\oplus r}$.
\qed

\neu{rank-gensur}{\bf Lemma. \/}{\em There exists a direct sum
$\bigoplus_{k=0}^{-\chi(E)}M_k$ of line bundles numerically equivalent to
$\Ocal_{\Pic^{0}(X)}$ and a commutative diagram with injective vertical
arrows}
\[ \xymatrix{\bigoplus_{k=0}^{-\chi(E)}M_k
\ar@{=}[r]\ar[d] & \bigoplus_{k=0}^{-\chi(E)}M_k \ar[d] \\
\FM^1(E(-D)) \ar[r]^\beta &
\FM^1(E) \,.} \]
\proof
Let $M$ be a line bundle on $X$ of degree $N$ for some $N \gg 0$.
By Serre duality and stability we conclude
$H^1(E(-D) \otimes M \otimes L)=0$ for
all line bundles $[L ] \in \Pic^0(X)$.
This is equivalent to $\FM^1_\Pcal(E(-D)\otimes M)=0$.
Taking a reduced section of $M$ and tensorizing with $E(-D)$
we obtain a short exact sequence
\[0 \to E(-D) \to E(-D) \otimes M \to \bigoplus_{i=1}^{N} k(P_i)^{\oplus r} \to 0.\]
Applying the Fourier-Mukai transform $\FM_\Pcal$ again, yields the exact
sequence
\[ \FM^0_\Pcal\left(\bigoplus_{i=1}^{N} k(P_i)^{\oplus r} \right)\to
\FM^1_\Pcal(E(-D)) \to 0 = \FM^1_\Pcal(E(-D)\otimes M) \,.\]
As seen in Lemma \ref{RANK}.\ref{rank-monad} the sheaf
$\FM^0_\Pcal\left(\bigoplus_{i=1}^{N} k(P_i)^{\oplus r} \right)$
is a direct sum of numerically trivial line bundles.
Again by Lemma \ref{RANK}.\ref{rank-monad} we obtain a sequence of surjections
\[ \FM^0_\Pcal\left(\bigoplus_{i=1}^{N} k(P_i)^{\oplus r} \right)\to
\FM^1_\Pcal(E(-D)) \to \FM^1_\Pcal(E) \,.\]
This implies the assertion because the rank of $\FM^1_\Pcal(E)$ is
at least $1-\chi(E)$ by assumption.
\qed

\neu{rank-h0s}{\bf Lemma. \/}{\em For a positive integer $m$ we have the
three equalities
\[ h^0\left(\bigoplus_{k=0}^{-\chi(E)} M_k(m \Theta) \right)=
(1-\chi(E))m^g \qquad h^0(\FM^0(\Ocal_D^{\oplus r})(m \Theta)) =
(g-1)r m^g \]
\[\mbox{and } \quad h^0(\FM_\Pcal^1(E(-D))(m \Theta)) =
 \left( g-1+\frac{g}{m} -\frac{\chi(E)}{r} \right) r m^g \,.     \]
}

\proof
The proofs of the statement for
$h^0\left(\bigoplus_{k=0}^{-\chi(E)} M_k(m \Theta) \right)$
and
$h^0(\FM^0(\Ocal_D^{\oplus r})(m \Theta))$
are the same by Lemma \ref{RANK}.\ref{rank-monad} and follow by basic facts
(see III.16 in Mumford's book \cite{Mum}) on the cohomology of line bundles
on abelian varieties.

Consider the vector bundle bundle
$F_m:=\Pcal \otimes q_X^*(E(-D)) \otimes p_X^* \Ocal(m \Theta)$
on $X \times \Pic^{0}(X)$.
The projection formula gives that $R^ip_{X*}F_m = \FM_\Pcal^i(E(-D))(m
\Theta)$.
Since $p_X$ has fiber dimension one, and $\FM_\Pcal^0(E(-D))=0$ by Lemma
\ref{RANK}.\ref{rank-monad} we obtain by the Leray spectral sequence
$H^0(F_m)=0 $ and the
isomorphism $H^1(F_m)  \cong H^0(\FM_\Pcal^1(E(-D))(m \Theta))$.

Next we apply the Leray spectral sequence to $q_X$ and the global section
functor.
It is shown in \cite{Muk} that $q_{X*}(\Pcal \otimes p_X^* \Ocal(m
\Theta))$ is the dual
of $P_m$. Thus, by the projection formula $q_{X*}(F_m) = P_m\dual \otimes
E(-D)$, and
$R^iq_{X*}(F_m)=0$ for all $i>0$.
We deduce from the spectral sequence that $H^0(F_m)= H^0(P_m\dual \otimes
E(-D))$ and
$H^1(F) \cong H^1(P_m\dual \otimes E(-D))$. Putting the result together we
find
\[ H^0(P_m\dual \otimes E(-D))=0 \qquad H^1(P_m\dual \otimes E(-D)) \cong
H^0(\FM_\Pcal^1(E(-D))) \,.\] 
Thus, $h^0(\FM_\Pcal^1(E(-D))) = h^1(P_m\dual \otimes E(-D)) = -
\chi(P_m\dual \otimes E(-D))$.
However, the last number can be directly computed by the Riemann-Roch
theorem for curves.
\qed

\neu{rank-h0}{\bf Corollary. \/}{\em
For $m> r \cdot g$ we have $h^0(\FM_\Pcal^0(E)(m \Theta)) > 0$.
}

\tip
\proof
From \ref{RANK}.\ref{rank-monad} and \ref{RANK}.\ref{rank-gensur} we
obtain 
the following diagram with exact row
\[\xymatrix{&&& H^0\left(\bigoplus_{k=0}^{-\chi(E)} M_k(m \Theta) \right)
\ar[d]_\gamma \\
0 \ar[r] & H^0(\FM_\Pcal^0(E)(m \Theta))  \ar[r] &
H^0(\FM_0(\Ocal_D^{\oplus r})(m \Theta)) \ar[r]^-{\alpha}
& H^0(\FM_\Pcal^1(E(-D))(m \Theta)) \ar[d]_\beta \\
&&& H^0(\FM^1_\Pcal(E)(m \Theta))\,. } \]
Furthermore, we have that $\beta \circ \gamma$ and $\gamma$ are injective
by Lemma \ref{RANK}.\ref{rank-gensur}, and
$\beta \circ \alpha = 0$ from the long exact sequence before
Lemma \ref{RANK}.\ref{rank-monad}.
Therefore the dimension of the image of $\alpha$ is at most
$h^0(\FM_\Pcal^1(E(-D))(m \Theta)) -
h^0 \left(\bigoplus_{k=0}^{-\chi(E)} M_k(m \Theta) \right)$.
This number is known by Lemma \ref{RANK}.\ref{rank-h0s} and by the
assumption on $m$ strictly
smaller than $h^0(\FM_0(\Ocal_D^{\oplus r})(m \Theta))$.
\qed

\neu{rank-unstable}{\bf Lemma. \/}{\em
If $E'$ is a vector bundle of with $\chi(E')>0$, and $rk(E') \leq r$,
then for all $m > r \cdot g$ we have $\Hom(P_m,E') \ne 0$.
}

\tip
\proof
We have that the slope of $E'$ is at least $\frac{1}{r}+(g-1)$.
The condition on $m$ implies that the slope of $E' \otimes P_m\dual$ is
strictly greater than $g-1$.
Thus, $\chi(E' \otimes P_m\dual) >0.$
This implies
$ 0 < \chi(E' \otimes P_m\dual) \leq h^0(E' \otimes P_m\dual) = \dim(\Hom(P_m,E'))$.
\qed

\neu{rank-thm2}{\bf Theorem. \/}{\em Let $E$ be a vector bundle on $X$
of rank $r$.
For any $m> r \cdot g$ we have the equivalence
\[ E \mbox{ has a } \Theta \mbox{-divisor} \quad  \Lrar \quad
\Hom(P_m,E)=0 \,.\]
}
\proof
If $E$ has a $\Theta$-divisor, then the claim follows from Theorem
\ref{PM}.\ref{pm-prop} and Lemma \ref{PM}.\ref{pm-chain}.

Let us assume that $E$ has no $\Theta$-divisor.
If for all $E' \subset E$ we have $\chi(E') \leq 0$, then we are done by
Corollary \ref{RANK}.\ref{rank-h0} and the identification of
$\Hom(P_m,E)$ with $H^0(\FM^0_\Pcal(E)(m \Theta))$ given in Theorem
\ref{PM}.\ref{pm-prop}.
If not, then $E$ contain a subsheaf $E'$ of positive Euler characteristic.
In this case we use the above Lemma \ref{RANK}.\ref{rank-unstable}.
\qed

\section{Minimal quotients of Raynaud's bundles}\label{BAD}\kkk{BAD}
\subsection{Minimal bundles}
\neu{bad1}{\bf $\chi$-small vector bundles.}
We have seen in Lemma \ref{RANK}.\ref{rank-unstable} that for a sheaf $E$
on $X$ of positive Euler characteristic $\chi(E)$
we have morphisms from $P_m$ to $E$ for $m > g \cdot \rk(E)$.
This implies, that for a vector bundle $E$ containing a vector bundle $E'$
of positive Euler characteristic
there are morphisms from $P_m$ to $E$ for $m \gg 0$.
Therefore, the existence of morphisms from $P_m$ to vector bundles which
do not contain any subsheaf of positive Euler characteristic
has to be investigated.
We call a sheaf $E$ which does not contain a subsheaf of positive 
Euler characteristic a $\chi$-small sheaf.
This implies that $E$ is a vector bundle. Semistable vector bundles of
slope $\leq g-1$ are
examples of $\chi$-small bundles. $\chi$-smallness can be expressed in terms
of the Harder-Narasimhan
filtration $0 \subset E_1 \subset E_2 \subset \ldots \subset E$ of $E$ by
the condition $\mu(E_1) \leq g-1$.

We frequently will use the obvious fact, that
subsheaves of $\chi$-small vector bundles are $\chi$-small too.

\neu{bad2}{\bf Morphisms from $P_m$ to $\chi$-small bundles.}
Let $X$ be a curve of genus $g \geq 2$. Since $P_2$ is $\chi$-small, we
have
morphisms from $P_m$ to $P_2$ for all $m \geq 2$ by Lemma
\ref{PM}.\ref{pm-chain}.
We consider the following set of integers 
\[N:=\left\{ \rk(E) \left| 
\begin{array}{l}
E \mbox{ is a $\chi$-small vector bundle,}\\
\Hom(P_m,E) \ne 0 \mbox{ for } m\gg 0\\
\end{array}
\right. \right\} \]
By the above remark $2^g = \rk(P_2) \in N$. Thus, $N \ne \emptyset$.
If $r \in N$, then there exists a $\chi$-small vector bundle $E$ of rank
$r$
with $\Hom(P_m,E) \ne 0$.
However, $E \oplus \Ocal_X$ is $\chi$-small too, and $\Hom(P_m,E)$ is
embedded into
$\Hom(P_m,E \oplus \Ocal_X)$. Eventually, we conclude $r+1 \in N$.
If $r_-$ is the minimal element of $N$, then we have seen, that 
$N=\{ r_- ,r_- +1,\ldots \}$. 
Of course, the number $r_-$ depends on the curve $X$.

By definition we have that a semistable vector bundle $E$ on $X$ with
$\chi(E)=0$ and $\rk(E) < r_-$ is not a base point for the generalized
$\Theta$-divisor.

\neu{bad3}{\bf Lemma.}
{\em If $E$ is a $\chi$-small vector bundle of rank $r_-$ on $X$
with $\Hom(P_m,E) \ne 0$ for $m\gg 0$, then there exists no non
trivial homomorphism
to a $\chi$-small vector bundle $E''$ of smaller rank.
Moreover, $E$ is stable.}\\
\proof
Suppose there exists a non trivial homomorphisms $E \rarpa{\alpha} E''$ to
a
$\chi$-small vector bundle $E''$ of rank $\rk(E'') < rk(E)=r_-$. Since the
image $\im(\alpha)$
is $\chi$-small too, we can assume that $\alpha$ is surjective. Setting
$E'=\ker(\alpha)$
we obtain another $\chi$-small vector bundle. Thus,
\[ 0 \to E' \to E \rarpa{\alpha} E'' \to 0\]
is a short exact sequence of $\chi$-small vector bundles.
The minimality of $r_-$ implies that $\Hom(P_m,E')=0=\Hom(P_m,E'')$. This
implies $\Hom(P_m,E)=0$
which contradicts our assumption.

Suppose $E$ is not semistable.
Take the first subsheaf $E_1 \subset E$ of the Harder-Narasimhan
filtration of $E$.
$E_1$ and $E/E_1$ are $\chi$-small of lower rank which is a contradiction.
If $E$ is semistable but not stable, then there exists a stable subbundle
$E_1 \subset E$
of the same slope. As before,
$E_1$ and $E/E_1$ are $\chi$-small. Thus, we conclude the stability of $E$.
\qed

\neu{bad4}{\bf Minimal bundles.}
We call a vector bundle $E$ on the curve $X$ a minimal bundle when
it satisfies the conditions

\tip
\begin{tabular}{ll}
(i) & $\Hom(P_m,E) \ne 0$ for $m \gg 0$, and\\
(ii) &  $\Hom(P_m,E') =0$ for all proper subsheaves $E' \subsetneqq E$.\\
\end{tabular}\\

By Theorems \ref{PM}.\ref{pm-prop} and \ref{RANK}.\ref{rank-thm2}
the following conditions are equivalent to
(i):

\begin{tabular}{ll}
(i-a) & $H^0(E \otimes L) \ne 0$ for all $[L ] \in \Pic^0(X)$,\\
(i-b) & $\Hom(L,E) \ne 0$ for all $[L ] \in \Pic^0(X)$,\\
(i-c) & $\Hom(P_m,E) \ne 0$ for all $m > g \rk(E)$.\\
(i-d) & $\Hom(P_m,E) \ne 0$ for one $m > g \rk(E)$.\\
\end{tabular}

\neu{bad4a}{\bf Lemma.}
{\em If $E$ is a sheaf on $X$ with $\Hom(P_m,E) \ne 0$ for $m \gg 0$,
then there exists a minimal subsheaf $F \subset E$.}\\
\proof We take an integer $m > g \rk(E)$. For all $\phi \in \Hom(P_m,E)$
the vector bundle $F_\phi:=\im(\phi)$ fulfills
$0 \leq \rk(F_\phi) \leq \rk(E)$,
and $\mu(P_m) \leq \mu(F_\phi) \leq \mu(E_1)$ where $E_1$ is the
sheaf from \ref{BAD}.\ref{bad2}. Consequently, the Hilbert polynomials of
the sheaves $\{ F_\phi \}_{\phi \in \Hom(P_m,E)}$ form a finite set.
Taking $F:=F_\phi$ with minimal possible Hilbert polynomial we obtain
the desired subsheaf by Theorem \ref{RANK}.\ref{rank-thm2}.
\qed

The subsheaf $F$ of the above lemma is not unique. However, to
understand sheaves $E$ with $\Hom(P_m,E) \ne 0$ it is convenient to
study the minimal ones. This is what we do next.

\subsection{The structure of minimal bundles}
\neu{bad-thm}{\bf Theorem.}
{\em Let $E$ be a minimal vector bundle on a smooth projective curve $X$
of genus $g$. Then there exists
an ample divisor $D=D(E)$ on $\Pic^0(X)$ such that
$\FM_\Pcal^0(i_*E) \cong \Ocal_{\Pic^0(X)}(-D)$, and
the vector bundle $P_D := \iota^* [-1]^*_{\Jac(X)} R^gq_* (\Pcal
\otimes p^* \Ocal_{\Pic^0(X)}(-D))$ admits a unique surjection to $E$.
}\\

\proof The proof of this theorem will be a consequence of the results in
lemma \ref{BAD}.\ref{bad6} -- \ref{BAD}.\ref{bad11}.
Throughout this subsection $E$ is a fixed minimal bundle.

\neu{bad5}{\bf Lemma.}
{\em Each nonzero morphism $\psi:P_m \to E$ to the minimal bundle $E$
is surjective.}\\
\proof
If $\psi:P_m \to E$ is not surjective, then the image of $\psi$ is a
vector bundle with a morphism from $P_m$ to it.
This contradicts the minimality of $E$.
\qed

\neu{bad6}{\bf Lemma.}
{\em For the minimal bundle $E$  the sheaf $\FM_\Pcal^0(E)$ is an ideal
sheaf $\Jcal_Z$.}\\
\proof
We consider a surjection $E \to k(x_0)$, and let $E'$ be the kernel.
We obtain an exact sequence
\[ 0 \to \FM_\Pcal^0(E') \to \FM_\Pcal^0(E) \to \FM_\Pcal^0(k(x_0)) = \Ocal_\Pic \,.\]
However, the minimality of $E$ implies that $\FM_\Pcal^0(E') = 0$, and 
$\FM_\Pcal^0(E) \ne 0$.
\qed

\neu{bad7}{\bf The divisor $D$.}
The ideal sheaf $\Jcal_Z$ can be decomposed as $\Ocal_\Pic(-D) \otimes \Jcal_{Z'}$
where $D$ is a divisor and $Z'$ is a subscheme of codimension greater than one.
This decomposition corresponds to the decomposition $Z = D \cup Z'$ of $Z$ into its
irreducible components of codimension one and those of greater codimension.

\neu{bad8}{\bf Lemma.}
{\em $D$ is an ample divisor.}\\
\proof
First $D$ is an effective divisor. Thus, the linear system $|2D|$ is base point free.
Thus, taking the Stein factorization of the morphism defined by the linear system
$|2D|$, we obtain a surjective morphism $\pi: \Pic^0(X) \to Y$ of projective
varieties with connected fibers. Let $T$ be the subscheme $\pi^{-1}(\pi(0))$ with
the reduced scheme structure. $D$ is ample if and only if $\pi$ is an isomorphism.
Thus, we have the obvious implication ($D$ is ample) $\Rar$
($T$ is a point). The converse implication holds too. Indeed,
assume that $T$ is a point. If all fibers of $\pi$ have dimension zero, then $D$ is
ample. Suppose that $C$ is curve in $\Pic^0(X)$ contained in a fiber of $\pi$. Thus,
$\deg(\Ocal(D)|_C) =0$. This holds true for all translates of the curve $C$. Let $C'$
be a translate of $C$ passing through $0$, then $C'$ must be contained in $T$ which
is a contradiction. Concluding, we have ($D$ is ample) $\Lrar$ ($T$ is a point).

$T$ is by definition a closed subscheme and closed under the group operation of
$\Pic^0(X)$. Hence, the embedding $\tau: T \to \Pic^0(X)$ is a morphism of abelian varieties.
We obtain a surjection $\tau^*: \Pic^0(\Pic^0(X)) \to \Pic^0(T)$ by
$[L] \mapsto \left[ (\Ocal(2D) \otimes L) |_T \right]$.

Let $[L] \in  \Pic^0(\Pic^0(X))$ be a line bundle on $\Pic^0(X)$ with $h^0(\Ocal(2D)
\otimes L) > 0$. Since, $\Ocal(2D) \otimes L$ is numerically trivial on the
fibers of $\pi$, the existence of a global section implies that it is trivial on the
generic fiber. Hence, by the seesaw theorem (see \cite{Mum}, p.~54) it is
trivial on all fibers, in particular it is trivial on $T$.
Thus, the set
$\{ [L] \in \Pic^0(\Pic^0(X)) \, | \, h^0(\Ocal(2D) \otimes L) > 0 \}$
is contained in the subgroup $\ker(\tau^*)$. 

Let $Q \in X$ be a geometric point, which maps under $\iota:X \to \Pic^0(\Pic^0(X))$
to the line bundle $L_Q$. As in the proof of lemma \ref{BAD}.\ref{bad6},
we conclude that there
exists an injection $\FM_\Pcal^0(E) \to L_Q$. Since $\Hom(\FM_\Pcal^0(E),L_Q) =
\Hom(\Ocal(-D),L_Q) = H^0(L_Q(D))$, we conclude that the line bundle $L_Q(D)$ has a
global section. Consequently, for any two points $Q,Q' \in X$ the line bundle $L_Q
\otimes L_{Q'}(2D)$ has a global section.

However, $\Pic^0(\Pic^0(X))$ is the smallest subgroup of $\Pic^0(\Pic^0(X))$ which
contains the image of $\iota:X \to \Pic^0(\Pic^0(X))$. Thus, the kernel of $\tau^*$
is  $\Pic^0(\Pic^0(X))$ itself. This implies that $\dim(T)=0$. Hence, the assertion of
the lemma holds.
\qed

\neu{bad9}{\bf Lemma.}
{\em The subscheme $Z'$ of \ref{BAD}.\ref{bad7} is empty. Hence, $\Jcal_Z=
\Ocal_{\Pic^0(X)}(-D)$.}\\
\proof 
We start with the short exact sequence on $\Pic^0(X)$:
\[ 0 \to \Jcal_Z \to \Ocal_{\Pic^0(X)}(-D) \to \Ocal_{\Pic^0(X)}(-D)|_{Z'} \to 0
\,,\]
and the resulting long exact sequence 
\[ \scriptstyle{
\dots \to \FM_\Pcal^{i-1} (\Ocal_{\Pic^0(X)}(-D)|_{Z'})
\to \FM_\Pcal^i(\Jcal_Z) \to \FM_\Pcal^i( \Ocal_{\Pic^0(X)}(-D)) \to
\FM_\Pcal^i (\Ocal_{\Pic^0(X)}(-D)|_{Z'})} \to \dots \]

Since the dimension of $Z'$ is at most $g-2$, we have
$\FM_\Pcal^i(\Ocal_{\Pic^0(X)}(-D)|_{Z'})=0$ for all $i> g-2$. The ampleness of $D$
(see lemma \ref{BAD}.\ref{bad8}) and Serre duality,
imply $\FM_\Pcal^i( \Ocal_{\Pic^0(X)}(-D))=0$, for all $i \ne g$.
Eventually, we conclude $\FM_\Pcal(\Ocal_{\Pic^0(X)}(-D)) = \FM^g_P(\Jcal_Z)[-g]$,
and the homomorphisms $\FM^g_P(\Jcal_Z) \to \FM_\Pcal^g(\Ocal_{\Pic^0(X)}(-D))$ is an
isomorphism.

The nontrivial morphism $\psi: \Jcal_Z \to \FM_\Pcal(\iota_*E)$ in $\Db(\Pic^0(X))$ 
induces a nonzero morphism
$\FM_\Pcal(\Jcal_Z)[g] \to [-1]^*\iota_*E=\FM_\Pcal(\FM_\Pcal(\iota_*E))[g]$.
The Eilenberg-Moore spectral sequence (see Theorem 2.11 in \cite{BK})
$E_2^{p,q}=\Ext^q(\FM_\Pcal^{g-p}(\Jcal_Z),[-1]^*\iota_*E) \Rar
\Ext^{p+q}(\FM_\Pcal(\Jcal_Z)[g],[-1]^*\iota_*E)$ is a first quadrant spectral sequence.
Thus, $\Hom(\FM_\Pcal(\Jcal_Z)[g],[-1]^*\iota_*E) =
\Hom(\FM^g_P(\Jcal_Z),[-1]^*\iota_*E)$.

Thus, we have in $\Db(\Jac(X))$ a commutative diagram
\[\xymatrix{ & \FM_\Pcal(\Ocal_{\Pic^0(X)}(-D)) \ar[dr]\\
\FM_\Pcal(\Jcal_Z) \ar[ur] \ar[rr]^-{\FM(\psi)} && [-1]^*\iota_*E[-g]\,.}\]
Since $[-1]^*\iota_*E[-g]=\FM_\Pcal(\FM_\Pcal(\iota_*E))$ and
$\FM_\Pcal$ is an equivalence,
we obtain a commutative diagram in $\Db(\Pic^0(X))$:
\[ \xymatrix{ & \Ocal_{\Pic^0(X)}(-D) \ar[dr]\\
\Jcal_Z \ar[ur] \ar[rr]^-\psi && \FM_\Pcal(\iota_*E) \,.} \] 
However, $H^0(\psi)$ is the identity of
$\Jcal_Z$. Since it factors though $\Ocal_{\Pic^0(X)}(-D)$,
we eventually yield the
stated equality $\Jcal_Z = \Ocal_{\Pic^0(X)}(-D)$.
\qed

\neu{bad10}{\bf The vector bundle $\Ecal$ on $\Jac(X)$.}
We define $\Ecal$ to be the line bundle
$[-1]^*\FM_\Pcal^g(\Ocal_{\Pic^0(X)}(-D))$.
From the morphism $\Ocal_{\Pic^0(X)}(-D) \to \FM_\Pcal(\iota_*E)$ we obtain
a nontrivial morphism $\xymatrix{\Ecal \ar[r]^-\pi & \iota_*E}$.

\neu{bad11}{\bf Lemma.}
{\em Properties of the vector bundle $\Ecal$.\\
\begin{tabular}{lp{14cm}}
(i) & The morphism $\pi:\Ecal \to \iota_*E$ is surjective.\\
(ii) & The rank of $\Ecal$ equals $h^0(\Ocal_{\Pic^0(X)}(D))$;\\
(iii) & If $[L] \in \Pic^0(\Jac(X))$ is the isomorphism class of a line bundle,
then $Hom(L, \Ecal)$ is of dimension one.\\
(iv) & If $\alpha:\Ecal \to k(Q)$ is a surjection to a skyscraper sheaf,
then those $[L] \in \Pic^0(\Jac(X))$ which admit a morphism to $\ker(\alpha)$
form a divisor in $\Pic^0(\Jac(X))$.\\
\end{tabular}
}
\proof
(i) Let $G \subset \iota_*E$ be the image of $\pi$. The sheaf $G$ is of the form
$\iota_*E'$ for a subsheaf $E' \subset E$. Since $\pi$ factors through $\iota_*E'$,
we have $\FM_\Pcal^0(\iota_*E')= \Ocal_{\Pic^0(X)}(-D)$. However, the minimality of
$E$ implies that for all proper subsheaves $F \subset E$, we have
$\FM_\Pcal^0(\iota_*F) =0$. Thus, $\pi$ is surjective.

(ii) Let $Q \in \Jac(X)$ be a closed point. As usual, $L_Q$ denotes the line bundle
on $\Pic^0(X)$ parameterized by $Q$. We have $\Hom(\Ecal,k(Q)) =
\Hom(\Ocal_{\Pic^0(X)}(-D),L_Q) = H^0(L_Q(D))$.

(iii) Since $L=\FM_\Pcal(k(R))$ for a point $R \in \Pic^0(X)$, we have
\[ \Hom(L,\Ecal) = \Hom(k(R),\Ocal_{\Pic^0(X)}(-D)[g]) = \Ext^g(k(R),
\Ocal_{\Pic^0(X)}(-D)) \,. \]
The sheaf $\Ocal_{\Pic^0(X)}(-D)$ is a line bundle. Therefore, we have
\[ \sext^i(k(R), \Ocal_{\Pic^0(X)}(-D)) \cong
\left\{ \begin{array}{ll} 0 & \mbox{ for } i \ne g\\
k(R) & \mbox{ for } i=g \, . \end{array} \right. \]
Therefore, the $H^k(\sext^l) \Rar \Ext^{k+l}$ spectral sequence degenerates, and we
conclude that $\Ext^g(k(R),\Ocal_{\Pic^0(X)}(-D))= H^0(k(R))$ is one dimensional.

(iv) As in $(ii)$ we have $\Hom(\Ecal,k(Q)) = \Hom(\Ocal_{\Pic^0(X)}(-D),L_Q)$.
Thus, for any non trivial morphism $\alpha: \Ecal \to k(Q)$, we obtain a nontrivial
morphism $\tilde \alpha: \Ocal(-D) \to L_Q$. Therefore $\tilde \alpha$ is injective.
The Fourier-Mukai transform of $\ker(\alpha)$ is given by the complex $\tilde \alpha:
\Ocal(-D) \to L_Q$. Now the statement follows from basic facts on cohomology and base
change (see II.5 in \cite{Mum}).
\qed

\neu{bad12}{\bf Lemma.}
{\em The surjection $\Ecal \to \iota_*E$ factors through $\Ecal|_C$. Let $F$ be the
kernel of $\Ecal|_C \to \iota_*E$. For the sheaf $F$ on $X$ we have
$\FM^0(\Ecal|_C) = \FM^0_\Pcal(F) \oplus \Ocal(-D)$.
}

\proof
(i) From the commutative diagram
\[\xymatrix{&&\Ecal \ar[d] \ar[dr]^\pi\\ 0 \ar[r] & F \ar[r] & \Ecal|_C \ar[r] &
\iota_*E \ar[r] & 0  } \]
and the fact that $\pi$ induces an isomorphism
$\FM^0_\Pcal(\pi): \FM^0_\Pcal(\Ecal) \to \FM^0_\Pcal( \iota_*E)$
it follows that
$\FM^0_\Pcal(\Ecal|_C) = \FM^0_\Pcal(F) \oplus \FM^0_\Pcal( \iota_*E)$.

\section{Base points for higher rank bundles}\label{SPEC}\kkk{SPEC}
\neu{spec-in} Let $X$ be a smooth projective curve of genus $g \geq 2$.
We say that a vector bundle $E$ with $\chi(E)=0$ has no $\Theta$-divisor
in $U_X(R,0)$ if for all vector bundles $F$ of rank $R$
the cohomology of $H^\bullet(X,E \otimes F)$ is not trivial.
By the Riemann-Roch theorem $\chi(E \otimes F)=0$ implies $\deg(F)=0$.
Furthermore from $H^\bullet(X,E \otimes F)=0$ we deduce the
semistability of $F$. Thus, denoting the set of all semistable vector
bundles of rank $R$ and degree zero on $X$ with $U_X(R,0)$ we see that
the above property is equivalent to 
\[ \Theta_{R,E} := \{ F \in U_X(R,0) \, | \, H^\bullet(X,E \otimes F)
\ne 0 \} = U_X(R,0) \,. \]

It is well known (see for example Beauville's survey articles \cite{Beau1}
and \cite{Beau2}) that: if $E$ is a base point of the
the $R$-th power of the generalized $\Theta$-divisor, then $E$ has no
$\Theta$-divisor in $U_X(R,0)$. Otherwise $\Theta_{R,E}$ is a divisor in
the moduli space $U_X(R,0)$.

The obvious fact that the direct sum of two semistable vector bundles
$F_1$ and $F_2$ of ranks $R_1$ and $R_2$ is semistable too yields:
If $E$ has a $\Theta$-divisor in $U_X(R_1,0)$ and in $U_X(R_2,0)$, then it has
a $\Theta$-divisor in $U_X(R_1+R_2,0)$.

\neu{spec-thm}{\bf Theorem. \/}{\em
There exists vector bundles $P_{m,R}$ on $X$ and canonical surjections
$P_{(m+1),R} \to P_{m,R}$ such that for a vector bundle $E$ of rank $r$
the conditions 
\begin{tabular}{ll}
(i) & $E$ has no $\Theta$-divisor in $U_X(R,0)$, \\
(ii) & $\Hom(P_{m,R},E) \ne 0$ for $m \gg 0$, and\\
(iii) & $\Hom(P_{m,R},E) \ne 0$ for $m> r (R^2(g-1)+1)$ 
\end{tabular}\\
are equivalent.
The numerical invariants of the bundle $P_{m,R}$ are given by
\[ \rk(P_{m,R}) = R m^{\tilde g} \quad
\deg(P_{m,R}) = \tilde g m^{\tilde g -1} \quad
\mu(P_{m,R}) = \frac{\tilde g}{mR} \quad
\mbox{ with }  \tilde g = R^2(g-1)+1\,. \]
}
\proof
Using the theory of spectral curves developed in \cite{BNR} we see that
there exists a finite morphism $f:\tilde X \to X$ of smooth projective
curves of degree $R$, such that\\
\begin{tabular}{lp{14cm}}
(i) & The direct image sheaf $f_*\Ocal_{\tilde X}$ is given by
$f_*\Ocal_{\tilde X} \cong \bigoplus_{k=0}^{R-1} \omega_X^{\otimes -k}$.
By Riemann-Roch, we have $\chi(\Ocal_{\tilde X})=R^2(1-g)$.
Thus, the genus $\tilde g$ of $\tilde X$ is given by
$\tilde g= R^2(g-1)+1$.\\
(ii) & Setting $\delta:=-\deg(f_*\Ocal_{\tilde X})=
\deg(\omega_{\tilde X} \otimes f^*\omega_X^{-1})=
(R^2-R)(g-1)$ we obtain a dominant morphism from an open subset $U
\subset \Pic^\delta(\tilde X)$  to $U_X(R,0)$ given by
$\tilde L \mapsto f_*\tilde L$.\\
\end{tabular}\\
This is Theorem 1 of \cite{BNR}. We note, that it also holds for
smooth projective curves defined over an algebraically closed field of
arbitrary characteristic.
Since the image of $f_*:U \to U_X(R,0)$ is not contained in a proper
closed subset, we conclude that $E$ has no $\Theta$-divisor in $U_X(R,0)$
if and only if for all $[\tilde L] \in U \subset
\Pic^\delta(\tilde X)$ we have $\Hom(f_* \tilde L,E) \ne 0$.
Since $[\tilde L] \mapsto \dim(\Hom(f_* \tilde L,E))$ is upper
semicontinuous we deduce that this is equivalent to
$\Hom(f_* \tilde L,E) \ne 0$ for all $\tilde L \in \Pic^\delta(\tilde
X)$. The functor $E \mapsto f^*E \otimes \omega_{\tilde X} \otimes
f^*\omega_X^{-1}$ is right adjoint to $\tilde L \mapsto f_* \tilde L$.
Thus, the above condition is equivalent to 
$\Hom(\tilde L,f^*E \otimes \omega_{\tilde X} \otimes f^* \omega_X^{-1})
\ne 0$ for all $\tilde L \in \Pic^\delta(\tilde X)$.
Fixing one $\tilde L_0 \in \Pic^\delta(X)$ we obtain 
that this is equivalent to $\Hom(\tilde M, \tilde L_0^{-1} \otimes
f^*E \otimes \omega_{\tilde X} \otimes f^* \omega_X^{-1}) \ne 0$ for all
$\tilde M \in \Pic^0(\tilde X)$. Now we can apply our Theorem
\ref{PM}.\ref{pm-prop} to deduce, that this is equivalent to
$\Hom(\tilde P_m, \tilde L_0^{-1} \otimes
f^*E \otimes \omega_{\tilde X} \otimes f^* \omega_X^{-1}) \ne 0$ for $m
\gg 0$ with $\tilde P_m$ the Raynaud bundle on $\tilde X$.
As before, we conclude that this is equivalent to 
$\Hom(f_*(\tilde P_m \otimes \tilde L_0),E ) \ne 0$ for $m \gg 0$.
Thus setting $P_{m,R}:=f_*(\tilde P_m \otimes \tilde L_0)$ we have
shown the equivalence of (i) and (ii).
Applying Theorem \ref{RANK}.\ref{rank-thm2} instead of Theorem
\ref{PM}.\ref{pm-prop} we obtain the equivalence of (i) and (iii).
The computation of the numerical invariants is straightforward.
\qed

\section{Applications and questions}\label{APP}\kkk{APP}
\subsection{Elliptic curves}
\neu{app-ell}{\bf Vector bundles on elliptic curves.}
Let $X$ be an elliptic curve. In this case we have (see Lemma
\ref{PM}.\ref{pm-num-inv}) the
numerical invariants
$\rk(P_m)=m$, and $\deg(P_m)=1$.
Therefore, the slope $\mu(P_m)=\frac{1}{m}$ is always positive.

Thus, for a vector bundle $E$ of rank $r$ and degree $0$ we deduce the
following equivalences
\[ E \mbox{ is semistable } \quad \Lrar \quad  \Hom(P_{r+1},E) =0 \quad
\Lrar \quad
H^*(X,E \otimes L)=0 \mbox{ for a } L \in \Pic^0(X).\]
Indeed, if $E$ is semistable, then from $\mu(P_{r+1})=\frac{1}{r+1}> 0 =
\mu(E)$ we derive that
$\Hom(P_{r+1},E)=0$. However, $\Hom(P_{r+1},E)=0$ implies by Theorem
\ref{RANK}.\ref{rank-thm2}
that $E$ has a $\Theta$-divisor. Having a $\Theta$-divisor implies semi
stability
immediately.

Moreover, Polishchuk shows in \cite{Pol} that the Fourier-Mukai transform
$\FM_\Pcal$
gives an equivalence between semistable bundles of rank $r$ and degree
zero
and torsion sheaves of length $r$.
See \cite{HP} for a presentation of Atiyah's results on vector bundles on
elliptic 
curves (cf.~\cite{Ati}) in terms of Fourier-Mukai transforms.

\subsection{Curves of genus two}
In this subsection we assume $X$ to be a curve of genus two defined over
the complex numbers. Recall that from \ref{BAD}.\ref{bad2} and
\ref{BAD}.\ref{bad3} that $r_-(X)$ is the smallest rank of a stable vector
bundle $E$ on $X$ with $\Hom(P_m,E) \ne 0$ for $m\gg 0$ and $\mu(E) \leq
g-1$.

\neu{app-g2-1}{\bf Theorem.\/}
{\em If $X$ is a curve $X$
of genus 2, then $r_-(X)=4$.}

\vspace{0.4em}
\proof
By  \ref{PM}.\ref{pm-num-inv} the Raynaud bundle $P_2$ is of rank four and has
$\mu(P_2)=1$. Thus, we have to exclude the ranks one, two, and three
as possible ranks of a stable vector bundle $E$ with $\Hom(P_m,E ) \ne
0$ and $\mu(E) \leq 1$.

{\bf Case 1: $\rk(E)=1$.\\}
This case can easily excluded because line bundles of degree
$d \leq g-1$ have a $\Theta$-divisor.

{\bf Case 2: $\rk(E)=2.$}\\
Here we have two subcases depending of the parity of the degree of $E$.

{\bf Case 2.1: $\rk(E)=2$ and $\deg(E)$ is even.}\\
Let $D$ be an effective divisor of degree $\frac{2g-2-\deg(E)}{2}$ and
take a global section $\xymatrix{\Ocal_X\ar[r]^-s & \Ocal_X(D)}$. 
This way we obtain a semistable vector bundle $E(D)= E \otimes \Ocal_X(D)$
with $\mu(E(D))=g-1$, and
the embedding $E \to E(D)$ gives $\Hom(P_m,E(D)) \ne 0$.
Thus $E(D)$ is a base point for the generalized $\Theta$-divisor which
is impossible for vector bundles of rank two as shown by Raynaud in
\cite{Ray} Corollaire 1.7.4.

{\bf Case 2.2: $\rk(E)=2$ and $\deg(E)$ is odd.}\\
Take an extension $0 \to E \to E' \to k(x_0) \to 0$ with $E'$ a vector
bundle. Since any subsheaf $L \subset E'$ gives a subsheaf $L(-x_0)$ of
$E$, the stability of $E$ implies the semistability of $E'$. Thus, we
pass to a semistable vector bundle $E'$ with $\deg(E')$ even and
$\Hom(P_m,E') \ne 0$. Now we proceed as in case 2.1.

{\bf Case 3: $\rk(E)=3.$}\\
Here we have three subcases according to $\deg(E)$ modulo 3.

{\bf Case 3.1: $\rk(E)=3$ and $\deg(E) \equiv 0 \mod 3$.}\\
Analogously to case 2.1 we obtain a base point for the generalized
$\Theta$-divisor on rank three bundles which contradicts again
Raynaud's result.

{\bf Case 3.2: $\rk(E)=3$ and $\deg(E) \equiv -1 \mod 3$.}
We proceed as in case 2.2.

{\bf Case 3.3: $\rk(E)=3$ and $\deg(E) \equiv 1 \mod 3$.}\\
Here we take a short exact sequence $0 \to E' \to E \to k(x_0) \to 0$.
The stability of $E$ implies as before the semistability of $E'$.
However $E$ is contained in $E'(x_0)$. The latter is semistable of
degree divisible by three. Hence case 3.1 applies.
\qed

\neu{app-g2-3}{\bf Corollary.\/}
{\em On a curve of genus 2 the Raynaud bundle $P_2$ is stable.}

\proof We have $\mu(P_2)=1$ by \ref{PM}.\ref{pm-num-inv}. Let $P_2 \to E$ be the
surjection to a stable bundle with $\mu(E) \leq 1$. By the above Theorem
\ref{APP}.\ref{app-g2-1} we have $\rk(E)=4$.
\qed

\subsection{Further examples}
\neu{app-g3}{\bf Curves of genus $g \geq 3$.}
Suppose that $X$ is a curve of genus 3 defined over the complex numbers.
Even though we can copy the proof of
Theorem \ref{APP}.\ref{app-g2-1} we only obtain $r_-(X) \geq 4$.
Here we use the base point freeness of the generalized $\Theta$-divisor
of bundles of rank two and three on the curve $X$ (see Proposition 1.6
in \cite{Beau2}).

Since base point freeness of the generalized $\Theta$-divisor on the
moduli space of rank three bundles is not known for a curve $X$ of genus
$g \geq 4$ we obtain a priori only $r_-(X) \geq 3$.

\neu{app-quot}{\bf Quot-schemes without torsion quotients.}
Let $X$ be a curve of positive genus $g$, and $m > g$ an integer.
Since any line bundle $L_0$ of Euler characteristic zero has a $\Theta$-divisor,
we conclude $\Hom(P_m,L_0)=0$. At the other hand we have $\Hom(P_m,L_1)
\ne 0$
for line bundles of positive Euler characteristic by Lemma
\ref{RANK}.\ref{rank-unstable}.

If $\psi:P_m \to L_1$ is a non-trivial morphism,
then it must be surjective. Otherwise it would factor through a subsheaf
$L_1(-D)$
for a non-trivial effective divisor $D$. This contradicts the fact, that
there exists no homomorphism from $P_m$ to a line bundle of Euler
characteristic zero. We
conclude
that the Quot-scheme $\Quot_X^{1,g}(P_m)$ of rank one quotients of
degree $g$ of $P_m$
has no torsion quotients. On the other hand, every line bundle $L_1$ of
degree $g$ appears
as a quotient of $P_m$ by Lemma \ref{RANK}.\ref{rank-unstable}.

\neu{app-quot2}{\bf Quot-schemes parameterizing only stable quotients.}
Let $X$ be a curve of genus $g \geq 2$.
The number $r_1$ is the smallest possible rank for a $\chi$-small vector
bundle $E$ with $\Hom(P_m,E) \ne 0$ for $m \gg0$.
Let $d_-$ be the smallest possible degree of such a bundle.

Then for $m \gg 0$ the Quot scheme
$\Quot^{r_-,d_-}_X(P_m)$ of quotients of $P_m$ of rank $r_-$ and degree
$d_-$
parameterizes only $\chi$-small vector bundles. Indeed, if $E$ is not
semistable, then it contains
a subsheaf $E_1 \subset E$ of maximal slope.
The surjection of $P_m \to E \to E/E_1$ gives us a surjection of $P_m$ to 
$E''=E/E_1$ with $rk(E'')<\rk(E)$ and $\deg(E'')< \deg(E)$. Proceeding
this way
we obtain a surjection from $P_m$ to a $\chi$-small bundle of degree less
than $r_-$
which is impossible by the very definition of $r_-$. Thus, each
quotient of $P_m$ with these
numerical invariants is a stable vector bundle. By definition of $r_-$ and
$d_-$
the scheme $\Quot^{r_-,d_-}_X(P_m)$ is not empty for $m \gg 0$.

\neu{app-mr}{\bf Generalized Raynaud bundles.}
For an ample line bundle $L$ on $\Pic^0(X)$ we can define
\[P_L := \iota^* [-1]^*_{\Jac(X)} R^gq_* (\Pcal \otimes p^* L^{-1})\,.\]
having the definition of the bundles $P_m$ in mind we obtain
$P_m=P_{m\cdot\Theta}$. Indeed, theorem \ref{BAD}.\ref{bad-thm}
implies that to each minimal bundle $E$ there exists a unique ample
divisor $D(E)$ on $\Pic^0(X)$ with a unique (up to scalars) surjection
$P_D \to E$. 

However, we gave another generalization of Raynaud's vector bundles with
the bundles $P_{m,R}$ of theorem \ref{SPEC}.\ref{spec-thm}.
To unify both, we consider a morphism $\pi: \tilde X \to X$
of irreducible smooth curves such that $\pi_*\Ocal_{\tilde X} \cong
\oplus_{k=0}^{R-1}\omega_X^{\otimes -k}$, and a line bundle $\tilde L_0$
on $\tilde X$ of degree $\delta=(R^2-R)(g-1)$.
Now any vector bundle $E$ with
\begin{itemize}
\item $\Hom(F,E) \ne 0$ for all vector bundles
$F$ of rank $R$ with $\chi(F)=0$; and
\item For all proper subsheaves $E' \subset E$, there exists a rank $R$
vector bundle $F$ with $\chi(F)=0$ and $\Hom(F,E')=0$.
\end{itemize}
determines an ample divisor $\tilde D$ in the Picard group
$\Pic^0(\tilde X)$ and a surjection
$\pi_*(P_{\tilde D} \otimes \tilde L_0) \to E$.
To proof this claim just combine the theorems \ref{BAD}.\ref{bad-thm}
and \ref{SPEC}.\ref{spec-thm}.

Having said this, it is natural to consider the vector bundles
$P_{\tilde D,R}:=\pi_*(P_{\tilde D} \otimes \tilde L_0)$ with $\tilde D$
an ample divisor on $\Pic^0(\tilde X)$ as the generalized Raynaud
bundles.

\neu{app-base}{\bf Base points of the generalized $\Theta$-divisor.}
To study the base points of the generalized $\Theta$-divisor on
$\U(r,r(g-1))$ it is by theorem \ref{RANK}.\ref{rank-thm2} enough to study
all quotients $Q$ of $P_{rg+1}$. Indeed, any base point corresponding to
the sheaf $E$ contains by this theorem a quotient $Q$ of $P_{rg+1}$ with
(semistability of $E$) Euler characteristic $\chi(Q) \leq 0$. The
sheaves $P_m$ with $1<m<gr+1$ are quotients of this type.
See \cite{Sch} for such a construction.

\neu{app-con}{\bf Question: Does $r_-(X)$ varies with $X$?}
It seems to me very probable that a curve $X$ where the N\'eron-Severi
group $\NS(\Pic^0(X))$ has high rank should have a smaller $r_-(X)$ than
a curve $X'$ with rank of $\NS(\Pic^0(X'))$ equal to one.

It is my hope that the generalized Raynaud bundles of
\ref{APP}.\ref{app-mr} will help to understand the dependence of
$r_-(X)$ from the N\'eron-Severi group $\NS(\Pic^0(X))$.

\vfill
{\small
Georg Hein, Universit\"at Duisburg-Essen, Fachbereich Mathematik,
45117 Essen\\
email: {\tt georg.hein@uni-due.de}}
\end{document}